\documentclass[12pt]{amsart}
\usepackage{amsmath,amssymb,amscd}

\theoremstyle{plain}
\newtheorem{theorem}{Theorem}

\newtheorem{lemma}[theorem]{Lemma}
\newtheorem{proposition}[theorem]{Proposition}
\newtheorem{definition}[theorem]{Definition}
\newtheorem*{conjecture}{Elliptic Stark Conjecture}

\theoremstyle{definition}
\newtheorem*{nonumtheorem}{Theorem}
\newtheorem*{remark}{Remark}

\newcommand{\thom}{\text{Hom}}
\newcommand{\tres}{\text{Res}}
\newcommand{\tind}{\text{Ind}}
\newcommand{\tdiv}{\text{div}}
\newcommand{\tdeg}{\text{deg }}
\newcommand{\tker}{\text{ker }}
\newcommand{\tdim}{\text{dim}\,}
\newcommand{\ttors}{{\text{tors}}}
\newcommand{\tord}{{\text{ord}}}
\newcommand{\tgal}{\text{Gal}}
\newcommand{\tnorm}{\text{Norm}}
\newcommand{\treg}{\text{reg}}

\begin{document}
\title[Stark conjectures for CM curves]{Stark conjectures \\for CM curves over number fields}
\author{Jeffrey Stopple}
\email{stopple@math.ucsb.edu}
\maketitle

\section{Introduction}\label{S:Intro}

In \cite{Birv}, Bloch constructs symbols in $K_2(E)$ for a CM elliptic curve $E$ defined over $\mathbb Q$,
corresponding to divisors  supported on torsion points of the curve.  This construction, and the special properties
of such curves, allowed him to prove the Beilinson conjecture for such curves.  In
\cite{D}, Deninger extends Bloch's results, for certain elliptic curves
 \lq of Shimura type \rq or \lq type (S) \rq .  For simplicity assume $E$ has complex multiplication by the ring of
integers
$\mathcal O_K$ of the complex quadratic field
$K$, and $E$ is defined over an extension $F$ of $K$. Shimura showed
\cite[Theorem 7.44]{Shim} that the following conditions are equivalent, and we will take either of them to mean
that
$E$ is of type (S).
\begin{nonumtheorem} $F(E_{\ttors} )$ is contained in
$K^{\text{ab}}$ if and only if  $F$ is abel\-ian over $K$ and the corresponding Hecke character $\psi$ on the
ideals of
 $F$ factors through the norm map from $F$ to 
$K$.
\end{nonumtheorem}

This condition is closely related to one considered by Gross in \cite{Gr}.
He calls a CM curve defined over a Galois extension $F$ of $\mathbb Q$ a 
\lq $\mathbb Q$-curve\rq if it is isogenous over $F$ to all
its Galois conjugates.  Similarly one defines $K$-curve via isogeny with all $\tgal(F/K)$ conjugates.
If $E$ is type (S) then it is a $K$-curve, since then the Hecke character $\psi$ is clearly Galois invariant and this is an isogeny
invariant \cite[Proposition (9.1.3)]{Gr}.  Conversely suppose $F$ is abelian over $K$ and $E$ is a $K$-curve defined over $F$.  If
$\mathcal O_K^\times=\{\pm1\}$ and the 2-Sylow subgroup of $\tgal(F/K)$ is cyclic, then $E$ is type (S) \cite[Proposition
2]{Robert}.

The key fact 
for results about special values of $L$-functions 
is that if $E$ is of type (S), then $L(s,E)$ factors as a product of $L$ functions of Hecke
characters of $K$.  In \cite{Ru}, results on the conjecture of Birch and Swinnerton-Dyer were obtained this way.  

To 
explain the Beilinson conjectures in the context of a curve 
$E$ defined over a number field $F$ in a simple format, we tensor 
$K_2(E_F)$ with $\mathbb Q$ to get a vector space $\mathbb Q K_2(E_F)$.  Let $\xi_i$ be a basis, and let $\Phi_j$ denote the
embeddings of $F$ into $\mathbb C$.  There are regulator maps for each embedding, described in \S \ref{S:Analysis} below
(along with some facts about $K$-theory).  Roughly speaking, the Beilinson conjectures  consists of
two parts.

\begin{enumerate}
\item\label{C:Dimension} \lq Dimension conjecture\rq: 
\[
\tdim \mathbb Q K_2(E_F)=[F:\mathbb Q]
\]
\item\label{C:Value} \lq $L$-value conjecture\rq:
\[
L(0,E)^{([F:\mathbb Q])}\approx_{\mathbb Q^\times}\det[\treg(\xi_i)_{\Phi_j}].
\]
\end{enumerate}
Note that for a CM curve,
the functional equation gives the order of the zero at $s=0$ as $[F:\mathbb Q]$.  What is known by the work of Bloch and
Deninger is that there are at least $[F:\mathbb Q]$ linearly independent symbols $\xi_i$ which make the $L$-value
conjecture hold.

In \S \ref{S:Torsion} below we prove a negative result: if the curve $E$ is not of Shimura
type, then the Bloch construction gives only symbols with regulator equal 0.   The problem is caused by the Galois action on the
torsion points.

In \S \ref{S:Stark} we try to use the Galois action on $K_2(E_F)$ to an advantage, by developing an elliptic curve analog of the
Stark conjecture.  
No claim is made that this conjecture is not
implied by more general motivic conjectures already in the literature.  The relevant $L$-functions are Artin-Hecke
$L$-functions.  Although we would like to work more generally, we restrict attention to CM curves for three reasons.  
First, for elliptic curves over number fields the continuation of the $L$-func\-tion to $s=0$ is still only conjectural.
More significant is that CM
curves have at worst additive bad reduction at any prime.  The $K$-theory of $E$ becomes more complicated if there is split
multiplicative bad reduction.   Finally,  a useful realization of the regulator map on $K_2(E)$
requires choosing a basis of the lattice in $\mathbb C$ corresponding to $E$.  For curves which are not defined over $\mathbb
R$ there is no canonical way to do this.  For CM curves we show in \S \ref{S:Analysis} that it is possible to make a choice so that
our determinant is well defined up to an element of $\mathbb Q^\times.$

In \S \ref{S:Rational}  assuming the Dimension conjecture (\ref{C:Dimension}) and the $L$-value conjecture
(\ref{C:Value}) above, we prove the analog
Stark's result
\cite{Stark} for rational characters.

In \S \ref{S:Abelian} assuming the Dimension conjecture and that $E$ is type (S), we prove the elliptic curve analog of the Stark
conjecture for an abelian extension of the complex quadratic field $K$.  In
particular, taking a trivial character of the Galois group, we have re-derived the result of Deninger in \cite{D} that
(\ref{C:Value}) above holds.  This is not an independent proof; the results rely on the same facts about curves of type (S) from
\cite{GS} that \cite{D} uses.  However, it is a very classical proof.  The two main ideas are an extension of the Frobenius
determinant theorem, and a distribution relation for values of Kronecker-Eisenstein series on isogenous curves.

In \S \ref{S:NotS} we prove a result for a general CM elliptic curve $E$ 
(i.e.; not necessarily type (S)) defined over an abelian extension $F$ of $K$:  There exists
an extension $M$ of $F$, and a character $\chi$ of $\tgal(M/F)$, such that the elliptic Stark conjecture is true for $E$ and
$\chi$.

The author would like to thank Fernando Rodriguez Villegas and Dinakar Ramakrishnan for helpful conversations, and Bill
Jacob for teaching the $K$-theory.

\section{Analytic prerequisites}\label{S:Analysis}

\subsection{The Bloch-Wigner Dilogarithm}
The regulator on $K_2$ of curves over number fields is  a map
into cohomology.
Here, however, we will follow the philosophy of \cite{Dinakar} where one finds the advice \lq\lq In
general, the more concrete one is able to make the [Borel] regulator map, the more explicit the information one is able to
extract from it.\rq\rq  So we will use Bloch's original, function theoretic approach to
the regulator as in \cite{Birv}.    Recall
that the classical Euler dilogarithm is defined by
\begin{align*}
Li_2(z)=&\sum_{n=1}^\infty \frac{z^n}{n^2}\qquad |z|<1\\
=&\int_0^z-\frac{\log(1-t)}{t} dt\qquad z\in\mathbb C\backslash[1,\infty)
\end{align*}
after analytic continuation.  The Bloch-Wigner dilogarithm
\[
D(z)=\text{Im}(Li_2(z))+\log|z|\arg(1-z)
\]
is well defined independent of path used to continue $Li_2$ and $\arg$.  For a torus
$\mathbb C/\Lambda$, $\Lambda=[\omega_1,\omega_2]$ corresponding to a point $\tau$ in $\mathcal H$,
we have the
$q$-symmetrized, or elliptic, dilogarithm
\[
D_q(z)=\sum_{n\in\mathbb Z} D(zq^n)\quad q=\exp(2\pi i\tau)\quad z\in \mathbb C^\times/q^{\mathbb
Z}\cong
\mathbb C/\Lambda.
\]
We define another real valued function $J(z)=\log|z|\log|1-z|$, and
\[
J_q(z)=\sum_{n=0}^\infty J(zq^n)-\sum_{n=1}^\infty
J(z^{-1}q^n)+\frac23\log|q|^2B_3(\log(\frac{|z|}{|q|})),
\]
where $B_3(t)$ is the third Bernoulli polynomial $B_3(t)=t^3-3/2t^2+1/2t$.
Together these functions make the regulator function: 
\[
R_q(z)=D_q(z)-iJ_q(z).
\]
(This normalization differs by $-i$ from the one usually taken.)
 
Recall that Weil \cite{W} defines the
Kronecker-Eisenstein-Lerch series for $x,x_0$ points in $\mathbb C/\Lambda$,
 and $\text{Re}(s)>a/2+1$ as
\[ 
K_a(x,x_0,s,\Lambda) = \sum\nolimits^*\langle x_0,\omega\rangle
(\overline{x}+\overline{w})^a |x+w|^{-2s}.
\] 
Here 
\[
\langle x_0,\omega\rangle=\exp((\overline{x_0}w-x_0\overline{w})/A(\Lambda)), 
\]
where 
$A(\Lambda)=\omega_1\overline{\omega_1}\text{Im}(\tau)/\pi$, so that
\[
\omega_1\overline{\omega_2}-\overline{\omega_1}\omega_2=-2\pi i\delta A(\Lambda)
\]
with $ A(\Lambda)>0$ and $\delta=\pm1$ chosen so that $\omega_2/\omega_1=\delta\tau$ with $\tau$ in $\mathcal H$.  In
this way $A(\Lambda)$ is independent of choice of generators for the lattice, while $A([1,\tau])=\text{Im}(\tau)/\pi$.
And $*$ means sum over $w\in\Lambda$, with $w\ne x$ if $x\in \Lambda$. 
We have the functional equation
\begin{multline}\label{E:fneqn}
\Gamma(s)K_a(x,x_0,s,\Lambda)=\\
A(\Lambda)^{a+1-2s}\Gamma(a+1-s)K_a(x_0,x,a+1-s)\langle x_0,x\rangle
\end{multline}
For the special case
$a=1$, $x=0$, $s=2$,  we will use the more concise notation
\[
K_{2,1}(u,\Lambda)=K_1(0,u,2,\Lambda)=\sum
\langle u,\omega\rangle\frac{\overline{w}}{|w|^4}.
\]
Observe that the behavior under homothety is simple:
\[
A^2(c\Lambda)K_{2,1}(cu,c\Lambda)=\overline{c}A^2(\Lambda)K_{2,1}(u,\Lambda).
\]

Bloch showed in
\cite{Birv} that for $\Lambda=[1,\tau]$
\begin{equation}\label{E:reg}
R_q(\exp(2\pi iu))=\pi A^2([1,\tau])K_{2,1}(u,[1,\tau]).
\end{equation}
It is worth observing how this function depends on the choice of the lattice basis.
Suppose we have a basis $\omega_1,\omega_2$ and $\tau=\omega_2/\omega_1$.  Let
\begin{gather*}
\tilde{\omega}_1=d\omega_1+c\omega_2\\
\tilde{\omega}_2=b\omega_1+a\omega_2
\end{gather*}
and $\tilde{\tau}=(a\tau+b)/(c\tau+d)$.  Then
\begin{align}\label{E:nwd}
\overline{\omega_1}A^2([1,\tau])K_{2,1}(u/\omega_1,[1,\tau])=&
A^2([\omega_1,\omega_2])K_{2,1}(u,[\omega_1,\omega_2])=\\
A^2([\tilde{\omega}_1,\tilde{\omega}_2])K_{2,1}(u,[\tilde{\omega}_1,\tilde{\omega}_2])=&
\overline{\tilde{\omega}_1}A^2([1,\tilde{\tau}])K_{2,1}(u/\tilde{\omega}_1,[1,\tilde{\tau}]).\notag
\end{align}
This is a problem in general, as there is no canonical choice for the lattice basis.

\subsection{Curves over number fields} \label{S:lambda}
Now let $E$ be an elliptic curve defined over a number field $F$, with complex multiplication by a
complex quadratic field $K$.   Let
$M$ be an extension of $F$, and $P$ a point in $E(M)$.  To each embedding
$
\Phi:M\hookrightarrow \mathbb C
$
we get an elliptic curve $E_\Phi$ over $\mathbb C$ corresponding to a lattice
$\Lambda$.  If $\Phi$ restricts to a real embedding of $F$, the real period gives a canonical choice for $\omega_1$. 
Extend to a lattice basis with any complex $\omega_2$ such that $\tau=\omega_2/\omega_1$ is in the upper half plane, let
$q=\exp(2\pi i\tau)$ and $u$ in $\mathbb C/\Lambda$ corresponding to $\Phi(P)$.  Then
\[
R_{q}(\exp(2\pi i u/\omega_1))=\pi A^2([1,\tau])K_{2,1}(u/\omega_1,[1,\tau])
\]
is well defined.
On the other hand, if $\Phi$ restricts to a complex embedding of $F$, we choose any basis $\omega_1,\omega_2$
for
$\Lambda$, and define $\tau$ and $q$ as before.  If $\Phi^\prime$ is the embedding which differs from $\Phi$
by complex conjugation, we can certainly choose basis
$\omega_1^\prime=\overline{\omega_1},\omega_2^\prime=-\overline{\omega_2}$ for
$\Lambda^\prime=\overline{\Lambda}$,  so $\tau^\prime=-\overline{\tau}$ and $q^\prime=\overline{q}$.

\begin{lemma}\label{relations} Suppose the embeddings $\Phi,\Phi^\prime$ differ by
complex conjugation.  Then
\begin{equation}\label{E:phiij}
R_{q}(\exp(2\pi iu/\omega_1))=\overline{R_{q^\prime}(\exp(2\pi iu^\prime/\omega_1^\prime))}
\end{equation}
\end{lemma}
\emph{Proof.} We have 
\begin{align*}
\text{Im}(\tau)^2K_{2,1}(u/\omega_1,\Lambda)=&
\text{Im}(\tau^\prime)^2
K_{2,1}(\overline{u^\prime}/\overline{\omega_1^\prime},\overline{\Lambda^\prime})\\
=&\text{Im}(\tau^\prime)^2\overline{K_{2,1}(u^\prime/\omega_1^\prime,\Lambda^\prime)}.
\end{align*}
So (\ref{E:phiij}) follows from (\ref{E:reg}).
\qed
\begin{remark}
Of course this still depends on the choice of the basis.  If we change $\tau$ to $(a\tau+b)/(c\tau+d)$, then (\ref{E:nwd})
implies that $R_{q}$ changes by $\overline{c\tau+d}$.
\end{remark}

To motivate what follows, we will summarize some relevant facts from $K$-theory.  For a commutative ring $R$, recall that
$K_0(R)$ is just the Grothendieck group, with generators $[M]$ for each projective $R$-module $M$, and relations
$[M]+[M^\prime]=[M\oplus M^\prime]$.  In particular for a field $k$, $K_0(k)=\mathbb Z$.  We will say no more about $K_1$
than the fact that for fields,
$K_1(k)=k^\times$.   For $K_2(k)$, the
Matsumoto relations give that
\[
K_2(k)=
k^\times\otimes k^\times/\{f\otimes 1-f\, |\,f\ne0,1\}.
\]
The class of $f\otimes g$, denoted $\{f,g\}$, is called a symbol.

For a curve $E$ over a field $F$, we have the function field $F(E)$, and a map
\begin{align*}
K_1(F(E))&\to\coprod_{P\in E(\overline{F})} K_0(\overline{F})=\text{Div}(E)\\
f&\mapsto \sum \tord_P f (P)
\end{align*}
The kernel $F^\times=K_1(F)$ of this map is relatively uninteresting; the cokernel $\text{Pic}(E)$ is important.

For  functions $f$, $g$ in $F(E)^\times$, and fixed $P\in E$, the tame symbol at $P$ is defined by
$$
T_P(f,g) = (-1)^{\tord_P g\,\tord_P f}\frac{f^{\tord_P g}}{g^{\tord_P f}}(P),
$$
and is trivial on tensors $f\otimes1-f$, thus is a function on symbols $\{f,g\}$.  In analogy to the divisor map above
we have
\[
K_2(F(E))\overset{\coprod T_P}{\longrightarrow}\coprod_{P\in E(\overline{F})} K_1(\overline{F})
\]
Here the cokernel is mysterious.  The kernel is, modulo torsion, our object of study $K_2(E)$.
 
When $F$ is a number field we get a group $K_2(E_\Phi)$ for each embedding $\Phi$ of $F$ into $\mathbb C$.
Associated to a symbol $\{f,g\}$ in $K_2(E_\Phi)$ we have the divisors $\tdiv(f)$, $\tdiv(g)$ of
the elliptic functions $f$, $g$.  Let 
\[ 
\sum_{Q,Q^\prime}\tord_Q(f)\tord_{Q^\prime}(g) (Q-Q^\prime)=\text{ say }\sum_Pa_P(P)
\]  
be their convolution 
$\tdiv(f)*\tdiv(g)$.  
Let $u_P$ be the point on $\mathbb C/\Lambda$ corresponding to $P$ on $E$.
The regulator associated to the symbol $\{f,g\}$ and the embedding
$\Phi$ is defined to be
\begin{align}
\treg(\{f,g\})_\Phi=&\sum_P a_PR_{q}(\exp(2\pi iu_{P}/\omega_1))/\pi\\
=&\sum_Pa_PA^2([1,\tau])K_{2,1}(u_{P}/\omega_1,[1,\tau]),\notag
\end{align}
where the dependence of each of the parameters $q$, $\tau$, $\omega_1$, and $u_P$ on the embedding $\Phi$ is suppressed. 
One can show \cite{Birv} that this is a Steinberg function, i.e. trivial on the relations that define $K_2$.

Now suppose that the number field $F$ has degree $n$ over $\mathbb Q$, and the extension $M$ is Galois over $F$ with Galois
group
$G$.   Let
$\Sigma=\thom(M,\mathbb C)$, and 
let $M_{\mathbb C}$ be the complex vector space with basis $\Sigma$.
A typical element in $M_{\mathbb C}$ is written $\sum_\Phi z_\Phi \Phi.$
Complex
conjugation acts on both
$\mathbb C$ and $\Sigma$, and we define Minkow\-ski space $M_{\mathbb R}$ to be the points such that
$z_{\overline{\Phi}}=\overline{z_\Phi}.$  This Euclidean space is canonically isomorphic to $\mathbb R^\Sigma$ \cite[chapter
I,\S5]{Neukirch}.  

We define a regulator map
\begin{gather*}
\lambda:K_2(E_M)\to M_{\mathbb R}\\
\xi\mapsto\sum\treg(\xi)_\Phi \Phi.
\end{gather*}
The relation (\ref{E:phiij}) of Lemma \ref{relations} show that $\lambda$ actually maps to $M_{\mathbb
R}$, not just
$M_{\mathbb C}$.
The action of $G$ on $\Sigma$ on the left is the opposite action on the field:
for $\gamma$ an element of $G$, 
$\gamma^{-1}\cdot\Phi(x)=\Phi(\gamma\cdot x)$, which we are writing $\Phi(x^\gamma)$.  So the $\Phi$ coefficient of
$\gamma\cdot\lambda(\xi)$, which is $\treg(\xi)_{\gamma^{-1}\cdot\Phi}$, is equal to $\treg(\xi^\gamma)_\phi$, the $\Phi$
coefficient of $\lambda(\xi^\gamma)$.  Thus the map $\lambda$ is a $G$ module homomorphism.

\begin{remark}
Of course, this $\lambda$ still depends on the choice of lattice basis at each embedding.   Suppose as before that
$\Phi,\Phi^\prime$ are related by complex conjugation and we have chosen $\tau$ and $\tau^\prime=-\overline{\tau}$. 
If we have a vector of symbols
$\overrightarrow{\xi}$, then by the remark after Lemma \ref{relations}, we compute that changing $\tau$ to
$(a\tau+b)/(c\tau+d)$ changes the vectors
\begin{align*}
\treg(\overrightarrow{\xi})_\Phi&\quad\text{ into }\quad
(c\overline{\tau}+d)\treg(\overrightarrow{\xi})_\Phi\\
\treg(\overrightarrow{\xi})_{\Phi^\prime}&\quad\text{ into }\quad
(c\tau+d)\treg(\overrightarrow{\xi})_{\Phi^\prime}.
\end{align*}
This changes the determinant of any matrix in which
these vectors appear, by $|c\tau+d|^2$.  Since our curve has complex multiplication, $\tau$ is in $K$ and this factor is in
$\mathbb Q^\times$.  Thus modulo $\mathbb Q^\times$, our determinants will be independent of choice of lattice basis.
\end{remark}

\subsection{Isogenies between curves}  Suppose we have elliptic curves $E$ and $E^\prime$ defined over a number field $F$,
and an $F$-isogeny 
\[
\phi:E\to E^\prime.  
\]
We identify the isogeny with a scalar
 $\phi\in\mathbb C^\times$ such that for the corresponding lattices,
$\phi \Lambda \subset \Lambda^\prime$ and
\begin{align*}
\phi:\mathbb C/\Lambda&\to\mathbb C/\Lambda^\prime\\
z&\mapsto \phi z
\end{align*}

The isogeny $\phi$ gives a contravariant map on the function fields of the two curves
\[
\phi^*f=f\circ\phi.
\]
This respects the Matsumoto relations so we get a map on $K_2$ of the function fields.  It is easy to show the tame symbol
satisfies
\[
T_P(\phi^*f,\phi^*g)=T_{\phi(P)}(f,g)
\]
so we get a map
\begin{align*}
\phi^*:K_2(E^\prime)&\to K_2(E)\\
\{f,g\}&\mapsto \{\phi^*f,\phi^*g\}.
\end{align*}

\begin{theorem}\label{P:distribution} We  have  the distribution relation
\[
\phi\cdot K_{2,1}(\phi(x),\Lambda^\prime)=\tdeg\phi\cdot
\sum_{t\in\text{ker}\phi}K_{2,1}(x-t,\Lambda).
\]
As a consequence, 
the regulator map $\lambda$ of \S\ref{S:lambda} satisfies
\[
\lambda(\phi^*\{f,g\})=\phi\cdot\lambda(\{f,g\}),
\]
 i.e. we have a commutative diagram
\[
\begin{CD}
K_2(E^\prime)   @> \phi^* >>  K_2(E)\\
@V \lambda VV   @VV \lambda V\\
F_{\mathbb R}  @> \phi >> F_{\mathbb R}
\end{CD}
\]
\end{theorem}
\emph{Proof}
We will first need a distribution relation for the isogeny given by multiplication by $d$ on $\mathbb C/\Lambda$ 
(found in \cite[Lemme 2.4.2]{MS}).
Fix $x$ in $\mathbb C/\Lambda$ and $x_0\in d^{-1}\Lambda$, then
\[
d^{2+a-2s}K_a(x_0,dx,s,\Lambda)=\sum_{t\in d^{-1}\Lambda/\Lambda}\langle-dx_0,x+t\rangle K_a(0,x+t,s,\Lambda).
\]
This is easy to prove as the left side is just
\[
d^2\sum\nolimits_{\omega_0}^* \langle \omega_0,dx\rangle \frac{(\overline{d\omega_0+dx_0})^a}{|d\omega_0+dx_0|^{2s}}
\]
while the right is just
\[
\sum_{\omega\ne-dx_0}\sum_t\langle \omega,x+t\rangle \frac{(\overline{\omega+dx_0})^a}{|\omega+dx_0|^{2s}}.
\]
The result follows from the orthogonality relation
\[
\sum_{t\in d^{-1}\Lambda/\Lambda} \langle \omega,t\rangle=
\begin{cases}
    d^2,&\text{ in case }\omega=d\omega_0\\
    0,&\text{ otherwise.}
 \end{cases}
\]
Now suppose $\Lambda,\Lambda^\prime$, and $\phi$ are as above, so $\phi\Lambda\subset\Lambda^\prime$.  By
the fundamental theorem of elementary divisors
(see, for example, \cite[Lemma 3.11]{Shim}),
we can choose bases $\{\omega_1,\omega_2\}$ of $\Lambda$ and
$\{\omega_1^\prime,\omega_2^\prime\}$ of $\Lambda^\prime$ so that
\[
\phi\omega_1=d_1\omega_1^\prime,\qquad\phi\omega_2=d_2\omega_2^\prime
\]
for some integers $d_1,d_2$, so 
\[
[\Lambda^\prime:\phi\Lambda]=d_1d_2=\text{deg}\phi.
\]
Note that $|\phi|^2A(\Lambda)=d_1d_2A(\Lambda^\prime)$, a fact we will use to compare the pairing
$\langle\, ,\, \rangle^\prime$ corresponding to $\Lambda^\prime$ with $\langle\, ,\,\rangle$.  We compute
\begin{align*}
K_a(0,\phi x,s,\Lambda^\prime)=&
\sum_{\omega^\prime\ne 0\in\Lambda^\prime}\langle \omega^\prime,\phi x\rangle^\prime
\frac{\overline{\omega^\prime}^a}{|\omega^\prime|^{2s}}\\
=&\frac{\overline{\phi}^a}{|\phi|^{2s}}
\sum_{\omega\ne 0\in\phi^{-1}\Lambda^\prime}\langle \omega, dx\rangle
\frac{\overline{\omega}^a}{|\omega|^{2s}}\\
=&\frac{\overline{\phi}^a}{|\phi|^{2s}}
\sum_{\tau\in\tker\phi}\langle\tau,dx\rangle K_a(\tau,dx,s,\Lambda)\\
=&\frac{\overline{\phi}^a}{|\phi|^{2s}}\tdeg\phi^{2s-a-2}
\sum_{\substack{\tau\in\tker\phi\\t\in\tker[\tdeg\phi]}}\langle [\tdeg\phi]\tau,t\rangle
K_a(0,x-t,s,\Lambda)
\end{align*}
by the distribution relation for the isogeny $[\tdeg\phi]$ above.  We use the orthogonality relation,
for each
$t\in\tker[\tdeg\phi]$
\[
\sum_{\tau\in \tker\phi} \langle [\tdeg\phi]\tau,t\rangle=
\begin{cases}
    \tdeg\phi,&\text{ if }t\in\tker\phi\\
    0,&\text{ if }t\in\tker[\tdeg\phi]\backslash \tker\phi
 \end{cases}
\]
to deduce
\[
K_a(0,\phi x,s,\Lambda^\prime)=\frac{\overline{\phi}^a}{|\phi|^{2s}}\tdeg\phi^{2s-a-1}
\sum_{t\in\text{ker}\phi}K_a(0,x-t,s,\Lambda).
\]
The distribution relation follows from taking $s=2$, $a=1$, and $\phi\overline{\phi}=\text{deg}\phi$.

For the commutative diagram, notice that if for each $Q$ we fix a $P$ in $\phi^{-1}(Q)$,
\[
\tdiv(\phi^*f)=\sum_Q\sum_{T\in\tker\phi}\tord_Q(f)(P-T),
\]
and similarly with $\tdiv(\phi^*g)$.
Thus
\[
\tdiv(\phi^*f)*\tdiv(\phi^*g)=\tdeg(\phi)\cdot \phi^*(\tdiv(f)*\tdiv(g)),
\]
and the result follows from the distribution relation.

\qed

\begin{remark}  Suppose the curve $E$ has complex multiplication by $\mathcal O_K$.  The theorem says the image of
$K_2(E)$ in $F_{\mathbb R}$ is an $\mathcal O_K$ module.
If the \lq Dimension conjecture\rq (\ref{C:Dimension}) is true,  this is a lattice in $F_{\mathbb R}$ with complex multiplication. 
\end{remark}

\section{Torsion divisorial support}\label{S:Torsion}

As mentioned in the introduction, a weak form of the Beilinson conjectures is known for CM elliptic curves of
type (S).  What can be done more generally?  Throughout this section we assume that $E$ is defined over an
abel\-ian extension $F$ of $K$, since the most interesting case is the \emph{smallest} extension of
$K$ over which $E$ is defined, the Hilbert class field $F=K(j(E))$.

Recall that $F(E_{\ttors} )$ is always an abel\-ian extension of $F$, and that
$K^{\text{ab}}$ is equal to $K(j(E),x(E_{\ttors}) ),$ where $x=x(P)$ is Weber's function for $E$.
Thus if $F$ is abel\-ian over $K$, so is $F(x(E_{\ttors}))$.

\begin{remark}  If the class number of $K$ is one, $K(j(E))=K$.  Then $F$ abel\-ian over $K$ implies
\[
F\subseteq K(x(E_{\ttors}) )\subseteq K(E_{\ttors} )
\]
so
\[
F(E_{\ttors} )\subseteq K(E_{\ttors} )\subseteq K^{\text{ab}}
\]
and so $E$ is of type (S).  We assume for the rest of this section that the class number of $K$ is
greater than one.\end{remark}

We would like to know when $\treg(\xi)_\Phi$
is equal to zero. In the paper \cite{S}, Schappacher considers this question for curves over
$\mathbb Q$.  Concerning the $d$ torsion, he remarks in (5.6)
that since the function $x\to K_1(x^\prime,x,s)$ is \emph{odd}, if a divisor $a$
is fixed under $a\to-a$, then  $K_1(0,a,s)=0$.  

Similar
consid\-erations ap\-pear in the thesis of Ross \cite{Ross}, from which we will borrow the idea of \lq
torsion divisorial support\rq:  Let $L$ be any extension of $F$, and let $i:F\hookrightarrow L$
denote inclusion.  This gives rise to two maps in $K$ theory, 
\[
i^*:K_2(E_F)\to K_2(E_L),
\]
 induced
by base extension, and 
\[
i_*:K_2(E_L)\to K_2(E_F),
\]
 induced by restriction of scalars.  The map
$i_*\circ i^*$ acts by multiplication by $[L:F]$, and $i^*\circ i_*$ is the norm map
\[
i^*\circ i_*:\{f,g\}\to \prod_{\sigma\in \tgal(L/F)}\{f,g\}^\sigma.
\]
If we fix an embedding $\Phi:L\hookrightarrow\mathbb C$, and let $\{f,g\}\in K_2(E_L)$, we see
\begin{equation}\label{E:istar}
\treg(i_*\{f,g\})_\Phi=\treg(i^*\circ i_*\{f,g\})_{\Phi|_F}=\sum_{\sigma\in
\tgal(L/F)}\treg(\{f,g\}^\sigma)_{\Phi|_F}.
\end{equation}
One sees immediately that the divisor corresponding to $i_*\{f,g\}$ is invariant under
$\tgal(L/F)$.  Ross then makes a definition similar to
\begin{definition}
Let 
$\mathcal N$ an ideal of $\mathcal O_K$.  Then a symbol in $K_2(E_F)$ is said
to have $\mathcal N$ torsion divisorial support if it is of the form
$i_*\prod_i\{f_i,g_i\}$
with all $f_i,g_i$ defined over $F(E_{\mathcal N})$, and 
such that the divisors $div(f_i),\, div(g_i)$ are supported on $E_{\mathcal
N}$.
\end{definition}
If there is a
$\sigma$ in $\tgal(F(E_{\mathcal N})/F)$ such that $P^\sigma=-P$ for all $\mathcal N$ torsion points
$P$, then we see the regulator is zero on any element $\xi$ with $\mathcal N$ torsion divisorial
support.  In this context the following
lemma will be useful.
\begin{lemma}\label{L:sc}  For an ideal $\mathcal N$ of $\mathcal O_K$ not dividing $2\mathcal O_K$,
the following are equivalent:
\begin{enumerate}
\item There exists an $\mathcal N$ torsion
point $P$ of $E$ such that
\[ 
\forall\sigma\in\tgal(F(E_{\mathcal N})/F),\quad P^\sigma\ne -P.\label{Sc}
\]
\item For all ideals $\mathcal A$ in
$\mathcal O_F$, 
\[
\psi(\mathcal A)\not\equiv -1 \mod \mathcal N, 
\]
where $\psi$ is the Hecke character of $E$
\label{Hecke}
\item $F(E_\mathcal N)=F(x(E_\mathcal
N)).$\label{Galois}
\end{enumerate}
\end{lemma}
\emph{Proof.}  If $\sigma$ is the Artin symbol of an ideal $\mathcal A$, then
$P^\sigma=\psi(\mathcal A)P$.  Thus \ref{Sc} and \ref{Hecke} are equivalent.
Since
$E$ is defined over
$F$, it is isomorphic over
$F$ to an equation of the form
$y^2=x^3+Ax+B$.  For a point $P=(x,y)$, we have $P^\sigma=-P$ if and only if $x^\sigma=x$ and
$y^\sigma=-y$.  If $F(E_\mathcal N)=F(x(E_\mathcal N))$ then clearly such a Galois action can not
happen.

Conversely suppose $P$ is $\mathcal N$ torsion such that for all $\sigma$, $P^\sigma\ne-P$, with
$\mathcal N$ minimal for $P$.  Let
$\sigma\in Gal(F(E_\mathcal N)/F(x(E_\mathcal N))).$  So for all $\tilde P$ in $E_\mathcal N$,
$x(\tilde P)^\sigma=x(\tilde P)$.  Choose a prime $\mathcal Q$ of $\mathcal O_F$ so that $\sigma$
is the Artin symbol for $\mathcal Q$; then $\tilde P^\sigma=\psi(\mathcal Q)\tilde P$.  From the
Weierstrass equation, if
$y(P)^\sigma\ne-y(P)$, it must equal $y(P)$.  So $\psi(\mathcal Q)P=P$, thus $\psi(\mathcal
Q)\equiv 1$ modulo $\mathcal N$.   Then for all $\tilde P$, $\psi(\mathcal Q)\tilde P=\tilde P$ and
$\tilde P^\sigma=\tilde P$, so $\sigma$ is trivial. This shows \ref{Sc} $\Leftrightarrow$
\ref{Galois}.
\qed

\begin{definition} $E$ is of type (R) if there exists an ideal $\mathcal N$ of $\mathcal O_K$ not
dividing $2\mathcal O_K$ such that any of the equivalent conditions above hold.
\end{definition} 

\begin{remark} This is a necessary condition for the regulator of a symbol with $\mathcal N$ torsion
divisorial support to be nonzero.
\end{remark} 

\begin{lemma}\label{L:ssc}  If $E$ is of type (S), then it is of type (R)\end{lemma}
\emph{Proof.}  This is Lemma 4.7 of \cite{GS}, where they show that (S) implies that for any
$\mathcal N$ divisible by both the conductor of $E$ and the conductor of $F$ over $K$, $F(E_\mathcal
N)=F(x(E_\mathcal N)),$ and is in fact the ray class field of $K$ modulo $\mathcal N.$
\qed

\begin{theorem}
$E$ is of type (R) if and only if it is of type (S).
\end{theorem}
\emph{Proof.}  Recall we are assuming the class number of $K$ is greater than 1, and thus $\mathcal
O_K^\times=\{\pm 1\}$.    By
\cite{Robert} Corollaire 2 we know there exists an elliptic curve $E^\prime$ defined over $F$ which
is of type (S).   The proof constructs a Hecke character $\psi^\prime$ which has the relevant
property.  By Theorem 9.1.3 of \cite{Gr}, we may assume that
$E^\prime$ and
$E$ have the same $j$ invariant.  Thus $E^\prime$ is a model of $E$ and so $\psi=\chi\psi^\prime$
for some quadratic Dirichlet character
$\chi$ associated to an extension $M/F$.  We will show that that $M$ is abel\-ian over $K.$

Let $\mathcal N$ be an ideal of $\mathcal O_K$ such that the Hecke character $\psi$ of $E$ is never
$-1$ modulo $\mathcal N$.  Since
$E^\prime$ is of type (S) it is of type (R) by Lemma \ref{L:ssc}, and we may assume there exists an
ideal $\mathcal N^\prime$ divisible by $\mathcal N$ such that $F(x({E^\prime}_{\mathcal
N^\prime}))$ is equal to
$F({E^\prime}_{\mathcal N^\prime})$.   Let $\mathcal Q$ a prime
ideal of $\mathcal O_F$ which splits completely in $F(x({E^\prime}_{\mathcal N^\prime}))$.  We will
show that
$\mathcal Q$ splits in $M$.  Since $Q$ splits completely, the corresponding Frobenius automorphism
$\sigma$ is trivial, so $P^\sigma=P$ for all $\mathcal N^\prime$ torsion $P$ on $E^\prime$, thus
$\psi^\prime(\mathcal Q)\equiv 1\mod \mathcal N^\prime$.  This means 
\[
\psi(\mathcal Q)=\chi(\mathcal
Q)\psi^\prime(\mathcal Q)\equiv \pm1\mod\mathcal N^\prime
\]
 because $\chi(\mathcal Q)=\pm1.$  Thus
$\psi(\mathcal Q)\equiv \pm1\mod\mathcal N$ as $\mathcal N$ divides $\mathcal N^\prime.$  By
hypothesis on
$\mathcal N$, we must have $\psi(\mathcal Q)\equiv 1\mod\mathcal N$, and therefore $\chi(\mathcal
Q)=1.$  Thus $\mathcal Q$ splits in $M$.  This then is enough to say that 
\[
M\subset
F(x({E^\prime}_{\mathcal N^\prime})),
\]
 a ray class field of $K.$  Thus $M$ is abel\-ian over $K$.  

By Lemme 1 of \cite {Robert}, we see that $E$ is of type (S).  The point is that
$E^\prime$ is of type (S) so  $F({E^\prime}_{\ttors} )$ is contained in
$K^{\text{ab}}.$  With $M$ also contained in $K^{\text{ab}},$ we get $F(E_{\ttors} )$ is
contained in $K^{\text{ab}}.$ \qed

\begin{remark} If the curve $E$ is not of Shimura type, then any symbol with torsion divisorial
support has regulator equal
$0$.
\end{remark}

\section{Stark conjectures}\label{S:Stark}

The result of the previous section presents two alternatives.  One could consider instead
the construction of symbols in $K_2(E)$ based on points of infinite order, as in \cite{GL}.    Or one can try to used the
Galois action to an advantage.  This possibility is suggested in \cite[p.187]{Dinakar}:
\lq\lq	 Because of the compatibility with the action of correspondences and Tate twists, the
Beilinson regulators admit a \lq motivic\rq\, formulation.  This generalizes Stark's conjectures on
the factoring of the regulator according to the Galois action, relating the eigen-pieces of unit
groups to the values of Artin $L$-functions at $s=0$...\rq\rq 

In this section we begin to work out the analog of Stark's conjectures for $K_2$ of an elliptic
curve with complex multiplication.   The approach here is as concrete and down to earth as
possible.  

\subsection{Notation}\label{SubS:Notation}  For any finite group $C$ and
class functions $\chi_1,\chi_2$ on $C$, we let $\langle\chi_1,\chi_2\rangle_C$ denote the scalar
product on $C$
\[
\langle\chi_1,\chi_2\rangle_C=\frac{1}{\sharp C}\sum_{t\in C} \chi_1(t)\overline{\chi_2(t)}.
\]
For any field $k$ and abel\-ian group $A$, we let $kA$ denote
$k\otimes A$.  Groups always act on the left, even if written $a^\sigma$ instead of $\sigma\cdot a$.  

Suppose $E$ is an elliptic curve over a number field $F$ with complex multiplication by $\mathcal
O_K$, where $K=\mathbb Q(\sqrt{D})$.  There are two cases:
\begin{enumerate}
\item $K\subseteq F$.  There exists a
Hecke character $\psi$ for $F$ such that
\[
L(s,E)=L(s,\psi)L(s,\overline{\psi}).\label{Case:Big}
\]
\item $K\nsubseteq F$; let $H=F\cdot K$.  
There exists a Hecke character $\psi$ for $H$ such that
\[
L(s,E)=L(s,\psi).\label{Case:Small}
\]
\end{enumerate}
We take a field $M$ Galois over $F$ which also contains $K$.  Let
$G=\tgal(M/F)$,
$N=\tgal(M/H)$, and $[F:\mathbb Q]=n$.

Fix an embedding $\Phi_1$ of $M$ (and thus also $H$, $K$) into $\mathbb C$, such that $\Phi_1(j(E))=j(\mathcal O_K)$.
Let $\Sigma_K=\thom_K(M,\mathbb C)$, so $\Sigma=\Sigma_K \cup \overline{\Sigma}_K.$
Recall the spaces $M_{\mathbb C}$, $M_{\mathbb R}$ of \S\ref{S:lambda}.

We need to define an appropriate $\mathbb Q$ vector space $M_{\mathbb Q}$ inside $M_{\mathbb R}$.  It seems
$M_{\mathbb Q}$ should be the field $M$ itself, viewed as $\sharp \Sigma_K$ copies of the field $K$.
The details of the $\mathbb Q[G]$ embedding inside $M_{\mathbb R}$ are in \ref{SubS:ChezStark} below.

The headings of the following subsections indicate which part of
\cite{Tate} we are imitating.

\subsection{$L$-Functions}

We take a finite
dimensional complex vector space $V$ and a representation
\[
\rho:G\to GL(V)
\] 
with character $\chi$.  Let $V^*$ denote the
contragredient representation of $V$. We define
\[
L(s,E\otimes\chi)=
	\begin{cases}
    L(s,\psi\otimes\chi)L(s,\overline{\psi}\otimes\chi),&\text{in case }\ref{Case:Big}\\
    L(s,\psi\otimes\tres_N\chi),&\text{in case }\ref{Case:Small}
 \end{cases}
\]
in terms of (products of) Artin-Hecke $L$ functions.   (Viewing the standard $L$ function of the elliptic curve as coming from a
Galois representation into cohomology, this is just the $L$ function of the tensor product representation.)
We see immediately this is well behaved with
respect to direct sums:
\begin{equation}
L(s,E\otimes(\chi_1\oplus\chi_2))=
L(s,E\otimes\chi_1)L(s,E\otimes\chi_2).
\end{equation}
For induction, we need the following fact about Artin-Hecke $L$ functions: if $F^\prime$ lies
between $F$ and $M$, fixed by $G^\prime$, and $\chi^\prime$ is the character of a representation of
$G^\prime$, then
\begin{equation}
L(s,\psi\otimes\tind_{G^\prime}^G(\chi^\prime))=
L(s,(\psi\circ\tnorm_{F^\prime/F})\otimes\chi^\prime).
\end{equation}
One shows then that
\begin{equation}
L(s,E_{F^\prime}\otimes\chi^\prime)=L(s,E_F\otimes\tind(\chi^\prime)).
\end{equation}
(There are three cases to check, depending on whether both $F^\prime$ and $F$ contain $K$, neither
do, or only $F^\prime$ does.)

\begin{proposition}  Let $n=[F:\mathbb Q]$.  Then
$L(s,E\otimes\chi)$ has a zero at $s=0$ of order $n\cdot\tdim(V)$.
\end{proposition}
\emph{Proof.}  Consider first case \ref{Case:Small}.  Via Brauer induction there exist one
dimensional characters $\chi_i$ on subgroups $G_i$ of $N$ and integers $n_i$ such that
\[
\tres_{N}(\chi)=\sum_i n_i \tind_{G_i}(\chi_i).
\]
Thus
\begin{align*}
L(s,\psi\otimes\tres_{N}(\chi))=&\prod_iL(s,\psi\otimes\tind_{G_i}(\chi_i))^{n_i}\\
=&\prod_iL(s,(\psi\circ\tnorm_{M^{G_i}/H})\otimes\chi_i)^{n_i}.
\end{align*}
Each of the $L$-functions
\[
L(s,(\psi\circ\tnorm_{M^{G_i}/H})\otimes\chi_i)
\]
has a zero at $s=0$ of order $[M^{G_i}:K]$ so the product has a zero of order 
\begin{align*}
\sum_i n_i[M^{G_i}:K]=&\sum_i n_i[M^{G_i}:H][H:K]\\
=&[H:K]\sum_i n_i \tdim(\tind_{G_i}(\chi_i))\\
=&[F:\mathbb Q]\tdim(V).
\end{align*}

In case \ref{Case:Big} similarly each of $L(s,\psi\otimes\chi)$ and
$L(s,\overline{\psi}\otimes\chi)$ have a zero at $s=0$ of order $\tdim(V)[F:K]$, so the
product has a zero of order $n\cdot\tdim(V)$.\qed

\begin{remark}
To get an appropriate regulator determinant, we want an automorphism of a vector space whose dimension is equal to the
order of the zero.   
\end{remark}
\begin{proposition} $n\cdot\text{dim}(V)=\tdim\,\thom_G(V^*,M_{\mathbb C})$
\end{proposition}
\emph{Proof.}
The representation of $G$
in $M_{\mathbb C}=M_{\mathbb Q}\otimes\mathbb C$ is just $n$ copies of the regular
representation
$\tind_e^G(1)$ of $G$, where $e$ denotes the identity element of $G$.
We have 
\begin{multline*}
n\cdot\text{dim}(V)=n\cdot\langle\tres_e(\chi),1\rangle_e=\\
\langle\tres_e(\chi),n\cdot1\rangle_e=
\langle\chi,n\cdot\tind_e^G(1)\rangle_G
\end{multline*}
by Frobenius Reciprocity.  Since the character of the right regular representation takes rational
integer values, it is real.   So the above is equal to
\[
\langle\chi\cdot n\cdot\tind_e^G(1),1\rangle_G=\tdim(V\otimes M_{\mathbb C})^G=
\tdim\,\thom_G(V^*,M_{\mathbb C})
\]
by duality.
\qed

\subsection{Stark Regulator}

\begin{proposition} \label{H:Hyp}  We suppose from now on 
the Dimension conjecture of \S \ref{S:Intro}; specifically,
that for all $L$,$F\subseteq L\subseteq M$, we have
\[
\tdim \mathbb Q K_2(E_L)=[L:\mathbb Q].
\]
Then $\mathbb Q K_2(E_M)\cong M_{\mathbb Q}$ as $\mathbb Q[G]$ modules.
\end{proposition}
\emph{Proof.}  Via the corollary on p.104 of \cite{Serre}, we need only show that for all subgroups
$C$ of $G$
\[
\tdim(\mathbb Q K_2(E_M)^C)=\tdim(M_{\mathbb Q}^C).
\]
Let $L$ the fixed field of $C$.  Via Galois descent for $K$ groups tensored with $\mathbb Q$,
\[
\mathbb Q K_2(E_M)^C\cong\mathbb Q K_2(E_L).
\]
By our assumption the dimension is $[L:\mathbb Q]$.  On the other hand,
\[
\tres_C\tind_e^G(1)=[G:C]\tind_e^C(1)
\]
by the Induction-Restriction theorem \cite[p.58]{Serre}.  This representation contains the trivial
representation $[G:C]=[L:F]$ times.  In $M_{\mathbb Q}$ we have $n=[F:\mathbb Q]$ copies of this
representation, so
\[
\tdim(M_{\mathbb Q}^C)=[L:\mathbb Q].
\]
\qed

\begin{remark}  One would like to try to get by with the weaker assumption 
$\tdim\,\mathbb Q K_2(E_L)\ge [L:\mathbb Q]$, since this is already known for curves of type
(S), by the work of Deninger \cite{D}.  One might hope to then prove that $M_{\mathbb Q}
\hookrightarrow \mathbb Q K_2(E_M)$ as $\mathbb Q[G]$ modules.  But the inequality on the
dimensions is not strong enough to prove this.  It might, for example, be true that
$\mathbb Q K_2(E_L)=\mathbb Q K_2(E_M)$ for every $L$, which would say
$\mathbb Q K_2(E_M)$ is trivial as a $\mathbb Q[G]$ module.
\end{remark}

Now, assuming the Dimension conjecture, let 
\[
f:  M_{\mathbb Q} \to \mathbb Q K_2(E_M)
\]
a $\mathbb Q[G]$ isomorphism.  
Recalling the map $\lambda$ defined in \S \ref{S:lambda}, we see
 $\lambda\circ f$ is a $G$ automorphism of $M_{\mathbb R}$, and
we use the same notation when extending scalars to $M_{\mathbb C}$.  Via functoriality, this defines an
automorphism $(\lambda\circ f)_V$:
\begin{gather*}
\thom_G(V^*,M_{\mathbb C})\overset{(\lambda\circ f)_V}{\to} \thom_G(V^*,M_{\mathbb C})\\
A \mapsto \lambda\circ f\circ A
\end{gather*}
We define $R(E,\chi)$ to be the determinant of $(\lambda\circ f)_V$.  (This determinant actually
depends on the choice of map $f$, which is suppressed from the notation.)  Let
$c(E,\chi)$ be the coefficient of the first term in the Taylor expansion of $L(s,E\otimes\chi)$
at $s=0$,
and define
\[
A(E,\chi)=\frac{R(E,\chi)}{c(E,\chi)}.
\]
\begin{conjecture} $A(E,\chi)$ belongs to $\mathbb Q(\chi)$, and for $\sigma$ in\newline $\tgal(\mathbb
Q(\chi)/\mathbb Q)$ we have
\[
 A(E,\chi)^\sigma=A(E,\chi^\sigma).
\]
\end{conjecture}

\subsection{Towards Stark's version}\label{SubS:ChezStark}

From our assumption on the dimensions of the $K$-groups in the previous subsection we deduced the
existence of a $\mathbb Q[G]$ isomorphism $f$ from $M_{\mathbb Q}$ to $\mathbb QK_2(E_M)$, which
implies there exists a set $\boldsymbol{M}$ of $n$ \lq Minkowski symbols\rq\,  
each of whose $G$ conjugates generate a $\sharp G$ dimensional $\mathbb Q$ vector space.
in which $G$ acts by the regular representation.  (These symbols
collectively play a role analogous to that of the Minkowski unit in the unit group.)  Conversely a set
$\boldsymbol{M}$ of
$n$ symbols in distinct
$G$ orbits will let us define a $\mathbb Q[G]$ isomorphism $f_{\boldsymbol{M}}$.  In this subsection we deduce a
more explicit formula for $R(E,\chi)$, by considering an explicit isomorphism $f_{\boldsymbol{M}}$.  

We must first specify the $\mathbb Q[G]$ embedding of $M$ inside $M_{\mathbb R}$.  Consider first case \ref{Case:Big}.  We
choose representatives for the $G$ orbits in
$\Sigma$ by taking
$\Phi_1$ as in \S \ref{SubS:Notation},
$\Phi_2,\dots,\Phi_{n/2}$ any representatives of the $G$ orbits in $\Sigma_K$, along with their complex conjugates.  Let $m$
in $M$ generate a normal basis of $M$ as an $F$ vector space, and $\{f_i\}$ any basis of $F$ as a $K$ vector space.  A typical
element of $M$ is then written
\[
\sum_{\sigma\in G}\sum_{i=1}^{n/2} c_i(\sigma)f_im^\sigma
\]
with all the $c_i(\sigma)$ in $K$, which we identify with
\[
\sum_{\sigma\in G}\sum_{i=1}^{n/2} c_i(\sigma)\sigma\cdot\Phi_i+\overline{c}_i(\sigma)\sigma\cdot\overline{\Phi}_i
\]
in $M_{\mathbb R}$.  To get a basis over $\mathbb Q$ we let 
\[
f_i^\pm=\frac{(1\pm\sqrt{D})}{2}f_i
\] 
for $i=1,\dots,n/2$.

Case \ref{Case:Small} is more complicated as $F$ is not a $K$ vector space.   
We write $n=r+2s$ with $r$ and $s$ the number of real embeddings, resp.
pairs of complex conjugate embeddings of $F$.  
We choose representatives 
$\Phi_1$ as above,
$\Phi_2,\dots,\Phi_r$
in $\Sigma_K$ so that $\Phi_i|F$ is a real
embedding for $i\le r$.  Thus $\overline{\Phi}_i$ is equal $\gamma_i^{-1}\cdot\Phi_i$ for some $\gamma_i$ in $G$.
For $r<i\le s$, $\Phi_i|F$ is a complex
embedding, so if we
fix once and for all a $\gamma$ in $G\backslash N$,
we can actually \textit{define} $\Phi_{i+s}=\gamma\cdot \overline{\Phi}_i$ in this case; then 
$\gamma^{-1}\cdot \overline{\Phi}_{i+s}=\Phi_i$.
Let $m$ in $M$ generate a normal basis of $M$ over $F$.   Let
$\{f_i\}$ a basis of $F$ as a $\mathbb Q$ vector space.  
Thus a typical element of $M$ can be written
\[
\sum_{\sigma\in G}\sum_{i=1}^n c_i(\sigma)f_im^\sigma
\]
with the $c_i(\sigma)$ in $\mathbb Q$.  For $i\le r$ we identify the element $f_im^\sigma$ in $M$ with
\[
(1+\sqrt{D})\,\sigma\cdot\Phi_i+(1-\sqrt{D})\,\sigma\cdot\overline{\Phi}_i.
\]
When $r<i\le s$, we identify the element $f_im^\sigma$ in $M$ with
\[
\sigma\cdot\Phi_i+\sigma\cdot\overline{\Phi}_i+
\sqrt{D}\,\sigma\cdot\Phi_{i+s}-\sqrt{D}\,\sigma\cdot\overline{\Phi}_{i+s}, 
\]
and $f_{i+s}m^\sigma$ with
\[
\sqrt{D}\,\sigma\cdot\Phi_i-\sqrt{D}\,\sigma\cdot\overline{\Phi}_i+
\sigma\cdot\Phi_{i+s}+\sigma\cdot\overline{\Phi}_{i+s}.
\]
We then extend $\mathbb Q$-linearly to get $M_{\mathbb Q}$ inside $M_{\mathbb R}$. 

The determinant of $(\lambda\circ f_{\boldsymbol{M}})_V$
in $\thom_G(V^*, M_{\mathbb C})$
 is equal to the determinant of
$1\otimes \lambda\circ f_{\boldsymbol{M}}$ in $(V\otimes M_{\mathbb C})^G$.
We're now ready to compute this determinant.

Consider first case \ref{Case:Big}.  We arbitrarily label the $n$ symbols as $\xi_{i,+}$
and $\xi_{i,-}$ for $i=1,\dots,n/2$.  We take the
isomorphism 
\[
f_{\boldsymbol{M}}(f_i^\pm m^\sigma)=\xi_{i,\pm}^\sigma\qquad i=1,\dots,n/2,\quad\sigma\in G.
\]
For $\tau\in G$ let  $\treg_{\boldsymbol{M}}(\tau)$ 
denote the $n\times n$ matrix of $2\times 2$ blocks
\[
\begin{bmatrix}
\treg(\xi_{i,+}^\tau)_{\Phi_j}&\treg(\xi_{i,+}^\tau)_{\overline{\Phi}_j}\\
\treg(\xi_{i,-}^\tau)_{\Phi_j}&\treg(\xi_{i,-}^\tau)_{\overline{\Phi}_j}
\end{bmatrix}
\] for $i,j=1,\dots,n/2$.   Then
\begin{proposition}\label{P:case1reg} In case \ref{Case:Big} we have
\begin{equation}\label{E:starkreg}
R(E,\chi,f_{\boldsymbol{M}})=(-\sqrt{D})^{\tdim(V)n/2}\det\left(\sum_{\tau\in
G}\rho(\tau)\otimes\treg_{\boldsymbol{M}}(\tau)\right).
\end{equation}
where as usual $\rho$ is the representation of $G$ with character $\chi$.
\end{proposition}
\emph{Proof.}
A typical vector in $(V\otimes M_{\mathbb C})^G$ looks like
\[
\sum_{\sigma\in G}\sum_{i=1}^{n/2}\rho(\sigma)v_i\otimes\sigma\cdot\Phi_i+
\rho(\sigma)v_i^\prime\otimes\sigma\cdot\overline{\Phi}_i
\]
with arbitrary vectors $v_i,v_i^\prime$ in $V$.  This space inherits a natural inner product from the one $\langle\, ,\,\rangle$
on $V$, namely the inner product of a typical vector as above with another, formed of vectors $w_i,w_i^\prime$ is just
\[
\sum_{i=1}^{n/2}\langle v_i,w_i\rangle+\langle v_i^\prime,w_i^\prime\rangle.
\]
We let $e_p$ for $p=1,\dots,\tdim(V)$ an orthonormal basis of $V$.
We choose a basis for $(V\otimes M_{\mathbb C})^G$ of the form
\[
v_{p,i}^\pm=\sum_{\sigma\in G} \rho(\sigma)e_p\otimes(f_i^\pm m^\sigma),
\]
for $p=1,\dots,\tdim(V)$, and $i=1,\dots,n/2$.  We are identifying $f_i^\pm m^\sigma$ with its image in $M_{\mathbb R}$ as in
\S \ref{SubS:Notation}.  This lets us compute 
\begin{align*}
1\otimes\lambda\circ f_{\boldsymbol{M}}(v_{p,i}^\pm)=&\sum_{\sigma\in G}\rho(\sigma)e_p\otimes
\lambda(\xi_{i,\pm}^\sigma)\\
=&\sum_{\sigma\in G}\rho(\sigma)e_p\otimes\sum_{\Phi\in\Sigma}\treg(\xi_{i,\pm})_\Phi \sigma\cdot\Phi.
\end{align*}
Write each $\Phi$ as $\tau^{-1}\cdot\Phi_k$ or $\tau^{-1}\cdot\overline{\Phi}_k$
and change variables $\sigma\mapsto\sigma\tau$ to get
\[
\sum_{\sigma\in G}\left\{\sum_{\tau\in G}\treg(\xi_{i,\pm}^\tau)_{\Phi_k}\rho(\tau)e_p\otimes\sigma\cdot\Phi_k+
\treg(\xi_{i,\pm}^\tau)_{\overline{\Phi}_k}\rho(\tau)e_p\otimes\sigma\cdot\overline{\Phi}_k\right\}.
\]
To get matrix coefficients we compute an inner product
\begin{multline*}
\langle 1\otimes\lambda\circ f_{\boldsymbol{M}}(v_{p,i}^\pm),v_{q,j}^\pm\rangle=
\langle \sum_{\tau\in G}\treg(\xi_{i,\pm}^\tau)_{\Phi_j}\rho(\tau)e_p,
\frac{(1\pm\sqrt{D})}{2}e_q\rangle\\
+\langle\sum_{\tau\in G}\treg(\xi_{i,\pm}^\tau)_{\overline{\Phi}_j}\rho(\tau)e_p,
\frac{(1\pm\sqrt{D})}{2}e_q\rangle
\end{multline*}
where the choices of $\pm$ on different sides of the inner product are of course independent.  Taking all four possible choices
of the $\pm$ gives us a $2\times 2$ block:
\[
\sum_{\tau\in G}\langle\rho(\tau) e_p,e_q\rangle
\begin{bmatrix}
\treg(\xi_{i,+}^\tau)_{\Phi_j}&\treg(\xi_{i,+}^\tau)_{\overline{\Phi}_j}\\
\treg(\xi_{i,-}^\tau)_{\Phi_j}&\treg(\xi_{i,-}^\tau)_{\overline{\Phi}_j}
\end{bmatrix}
\begin{bmatrix}
\frac{1-\sqrt{D}}{2}&&\frac{1+\sqrt{D}}{2}\\
\frac{1+\sqrt{D}}{2}&&\frac{1-\sqrt{D}}{2}
\end{bmatrix}
\]

The determinant of the matrix with these (doubly indexed) coefficients is our regulator
$R(E,\chi,f_{\boldsymbol{M}})$.
\qed

Case \ref{Case:Small} seems, at first, simpler.  We have symbols $\xi_i$
 for $i=1,\dots,n$ in distinct $G$ orbits.  We take the
isomorphism 
\[
f_{\boldsymbol{M}}(f_i m^\sigma)=\xi_i^\sigma\qquad i=1,\dots,n,\quad\sigma\in G.
\]
Now a typical vector in $(V\otimes M_{\mathbb C})^G$ looks like
\[
\sum_{\sigma\in G}\sum_{i=1}^{n}\rho(\sigma)v_i\otimes\sigma\cdot\Phi_i
\]
with arbitrary vectors $v_i$ in $V$.  The inner product of a typical vector as above with another formed of vectors
$w_i$ is just
\[
\sum_{i=1}^{n}\langle v_i,w_i\rangle.
\]
We let $e_p$ for $p=1,\dots,\tdim(V)$ an orthonormal basis of $V$.
We choose a basis for $(V\otimes M_{\mathbb C})^G$ of the form
\[
v_{p,i}=\sum_{\sigma\in G} \rho(\sigma)e_p\otimes(f_i m^\sigma),
\]
for $p=1,\dots,\tdim(V)$, and $i=1,\dots,n$.  This lets us compute 
\begin{align*}
1\otimes\lambda\circ f_{\boldsymbol{M}}(v_{p,i})=&\sum_{\sigma\in G}\rho(\sigma)e_p\otimes
\lambda(\xi_i^\sigma)\\
=&\sum_{\sigma\in G}\rho(\sigma)e_p\otimes\sum_{\Phi\in\Sigma}\treg(\xi_i)_\Phi \sigma\cdot\Phi\\
=&\sum_{\sigma\in G}\left\{\sum_{\tau\in G}\treg(\xi_{i}^\tau)_{\Phi_k}\rho(\tau)e_p\right\}\otimes\sigma\cdot\Phi_k
\end{align*}
after writing each $\Phi$ as $\tau^{-1}\cdot\Phi_k$ 
and changing variables $\sigma\mapsto\sigma\tau$ just as before.

In order to compute matrix coefficients 
as inner products, we must rewrite the typical basis vector $v_{q,j}$ with $j\le r$ as
\begin{align*}
v_{q,j}&=\sum_{\sigma\in G}\rho(\sigma) e_q\otimes((1+\sqrt{D})\sigma\cdot\Phi_j+
(1-\sqrt{D})\sigma\cdot\overline{\Phi_j})\\
&=\sum_{\sigma\in G}\rho(\sigma) e_q\otimes((1+\sqrt{D})\sigma\cdot\Phi_j+
(1-\sqrt{D})\sigma\gamma_j^{-1}\cdot\Phi_j)\\
&=\sum_{\sigma\in G}\rho(\sigma)((1+\sqrt{D})e_q+(1-\sqrt{D})\rho(\gamma_j)e_q)\otimes\sigma\cdot\Phi_j
\end{align*}
after using the relation $\overline{\Phi}_j=\gamma_j^{-1}\Phi_j$ of \ref{SubS:Notation} above, and a change of variables.  We
then see that for $j\le r$
\begin{multline*}
\langle 1\otimes\lambda\circ f_{\boldsymbol{M}}(v_{p,i}),v_{q,j}\rangle=\\
\sum_{\tau\in G}\treg(\xi_i^\tau)_{\Phi_j}\langle\rho(\tau)e_p,(1+\sqrt{D})e_q+(1-\sqrt{D})\rho(\gamma_j)e_q\rangle\\
=(1-\sqrt{D})\sum_{\tau\in G}\treg(\xi_i^\tau)_{\Phi_j}\langle\rho(\tau)e_p,e_q\rangle+\\
(1+\sqrt{D})\sum_{\tau\in G}\treg(\xi_i^\tau)_{\Phi_j}\langle\rho(\tau)e_p,\rho(\gamma_j)e_q\rangle
\end{multline*}
In the second sum above we use that $\rho$ acts by isometries
\[
\langle\rho(\tau)e_p,\rho(\gamma_j)e_q\rangle=\langle\rho(\gamma_j^{-1}\tau)e_p,e_q\rangle
\]
and change the variables $\tau\mapsto \gamma_j\tau$.  Then
\[
\treg(\xi_i^{\gamma_j\tau})_{\Phi_j}=\treg(\xi_i^\tau)_{\gamma_j^{-1}\cdot\Phi_j}=
\treg(\xi_i^\tau)_{\overline{\Phi}_j}
\]
so finally
\begin{multline*}
\langle 1\otimes\lambda\circ f_{\boldsymbol{M}}(v_{p,i}),v_{q,j}\rangle=\\
\sum_{\tau\in G}
\left\{(1-\sqrt{D})\treg(\xi_i^\tau)_{\Phi_j}+(1+\sqrt{D})\treg(\xi_i^\tau)_{\overline{\Phi}_j}\right\}
\langle\rho(\tau)e_p,e_q\rangle
\end{multline*}
for $j\le r$.

Similarly if $r<j\le r+s$, we rewrite $v_{q,j}$ as
\begin{multline*}
v_{q,j}=\sum_{\sigma\in G}\rho(\sigma)(e_q-\sqrt{D}\rho(\gamma^{-1})e_q)\otimes\sigma\cdot\Phi_j+\\
\rho(\sigma)(\sqrt{D}e_q+\rho(\gamma)e_q)\otimes\sigma\cdot\Phi_{j+s})
\end{multline*}
using the relation between $\Phi_j$ and $\overline{\Phi}_{j+s}$ above and a change of variables.
Using the same change of variables as above, we see that
\begin{multline*}
\langle 1\otimes\lambda\circ f_{\boldsymbol{M}}(v_{p,i}),v_{q,j}\rangle=
\sum_{\tau\in G}\langle\rho(\tau)e_p,e_q\rangle\times\\
(\treg(\xi_i^\tau)_{\Phi_j}+\treg(\xi_i^\tau)_{\overline{\Phi}_j}
-\sqrt{D}\treg(\xi_i^\tau)_{\Phi_{j_s}}+\sqrt{D}\treg(\xi_i^\tau)_{\overline{\Phi}_{j+s}}).
\end{multline*}

\section{rational characters}\label{S:Rational}
For the regulator determinant we clearly we have
\[
R(E,\chi_1\oplus\chi_2)=R(E,\chi_1)R(E,\chi_2).
\]
Induction properties, as well, seem to follow from \cite[p.29]{Tate}: If $\chi$ is the character of a representation of
$\tgal(M/F^\prime)$, then
\[
R(E_{F^\prime},\chi)=R(E_F,\tind(\chi)).
\]
For example, with $F^\prime=M$ we get
\[
R(E_M)=R(E_F,\tind_e^G(1))=\prod_{\chi\in\hat{G}}R(E_F,\chi)^{\tdim(\chi)},
\]
So we have factored the regulator determinant into pieces.   By the usual properties of $L$-functions, the same holds for
the $c(E,\chi)$ and thus for the ratios $A(E,\chi)$.

\begin{remark}  If we take the trivial representation $1$ of $G$,  then
comparing (\ref{E:istar}) and (\ref{E:starkreg}) we see we have recovered the map $i_*:K_2(E_M)\to
K_2(E_F)$, and the determinant $R(E_F)$ is one piece of the determinant $R(E_M)$.  As in \S \ref{S:Torsion},
if $E$ is not of type (S), the regulator map is zero on symbols with torsion divisorial support.  The  other terms
in the product are not, however, \emph{a priori} zero on symbols with torsion divisorial support.  For
example, take $F=K(j(E))$ the Hilbert class field, and
$M=F(E_{\mathcal Q})$, where $p\mathcal O_K=\mathcal Q\overline{\mathcal Q}$ is a split prime, and
take $\rho$ an odd character of
$\tgal(M/F)=\mathbb F_p^\times$.  This is a simple observation, but, as mentioned at the beginning of this
section, it is the motivation for looking at this analog of Stark's conjecture.
\end{remark}

\begin{theorem}  Suppose $E$ is an elliptic curve defined over $F$, and $M$ is a Galois
extension, with $\chi$ a character  of a representation of $\tgal(M/F)$ taking rational values.   Then there
exists integers $m,n_i$ and intermediate fields $F_i$ such that
\[
L(s,E_F\otimes\chi)^m=\prod_i L(s,E_{F_i})^{n_i}.
\]
If we assume the Dimension conjecture and the $L$-value conjecture of \S \ref{S:Intro}, we get that
\[
A(E,\chi)^m=\prod_i A(E_{F_i})^{n_i}\in\mathbb Q.
\]
\end{theorem}
\emph{Proof.}  
We will treat the case when $F$ does not
contain $K$, the other case is easier.  Let $m$ be the exponent of $G=\tgal(M/F)$.  By standard facts about
representations (e.g.
\cite[p.103]{Serre})  there are subgroups $C_i$ of $G$ and integers $n_i$ such that
\[
m\cdot\chi=\sum_i n_i\tind_{C_i}^G(1),\qquad 
m\cdot\tres_N(\chi)=\sum_i n_i\tres_N\tind_{C_i}^G(1).
\]
So
\[
L(s,E\otimes m\cdot\chi)=\prod_i L(s,\psi\otimes\tres_N\tind_{C_i}^G(1))^{n_i}.
\]
We want to use the Induction-Restriction theorem \cite[p.58]{Serre} on each term, so we need for each $i$ a
decomposition of 
$G$ into double cosets
\[
G=\cup_\gamma N\gamma C_i.
\]
Fix an element $\delta\in G$, $\delta\notin N$;
since $[G:N]=2$, there are at most two double cosets $Ne C_i$ and $N\delta C_i$.  There are two cases:  
\begin{enumerate}
\item $C_i\nless N$;  Then there is only a single double coset $Ne C_i=N\delta C_i$.  Let
$\tilde{C}_i=C_i\cap N$, $F_i$ the fixed field of $C_i$ and $\tilde{F}_i$ the fixed field of $\tilde{C}_i$.  The
Induction-Restriction theorem says
\[
\tres_N\tind_{C_i}^G(1)=\tind_{\tilde{C}_i}^N(1).
\]
So 
\begin{align*}
L(s,\psi\otimes\tres_N\tind_{C_i}^G(1))=&
L(s,\psi\otimes\tind_{\tilde{C}_i}^N(1))\\=&
L(s,\psi\circ\tnorm^{\tilde{F}_i}_H)\\=&L(s,E_{F_i}).
\end{align*}
\item $C_i<N$.  The two double cosets are distinct.  Let $D_i=\delta C_i\delta^{-1}$, also a subgroup of $N$.
Let $F_i$ the fixed field of $C_i$, and $L_i$ the fixed field of $D_i$.
The
Induction-Restriction theorem says
\[
\tres_N\tind_{C_i}^G(1)=\tind_{C_i}^N(1)\oplus \tind_{D_i}^N(1).
\]
So 
\begin{align*}
L(s,\psi\otimes\tres_N\tind_{C_i}^G(1))=&
L(s,\psi\otimes\tind_{\tilde{C}_i}^N(1))L(s,\psi\otimes\tind_{\tilde{D}_i}^N(1))\\=&
L(s,\psi\circ\tnorm^{F_i}_H)L(s,\psi\circ\tnorm^{L_i}_H)\\
\intertext{Since $E$ is defined over the subfield $F$ of $H$, and $\delta$ generates $\tgal(H/F)$, we get
from
\cite[Theorem 10.1.3]{Gr} that
$\delta\cdot\psi=\overline{\psi}$.  This gives}
=&L(s,\psi\circ\tnorm^{F_i}_H)L(s,\overline{\psi}\circ\tnorm^{F_i}_H)\\
=&L(s,E_{F_i}).
\end{align*}
\end{enumerate}

Assuming the $L$-value conjecture of \S \ref{S:Intro}, there exist rational numbers
\[
A(E_{F_i})=\frac{R(E_{F_i})}{c(E_{F_i})}.
\]
Then the induction and direct sum properties imply the following weak form of the elliptic Stark conjecture for rational
characters:
\[
A(E,\chi)^m=\prod_i A(E_{F_i})^{n_i}\in\mathbb Q.
\]
\qed

\section{abel\-ian over complex quadratic}\label{S:Abelian}

This section is devoted to the proof of the following
\begin{theorem}
Suppose the field $M$ has abelian Galois group $\Gamma$ over the complex quadratic
field $K$,  and that the curve $E$ is type (S).  Then the elliptic Stark conjecture is true.   
\end{theorem}
To minimize notation, we will consider the case where $E$ is
defined over $F=K(j(E))$.  This makes $F$ the Hilbert class field of $K$,
and $n=[F:\mathbb Q]=2h$ where $h$ is the class number of $\mathcal O_K$.  
For abelian Galois groups we may as well assume the character $\chi$ satisfies $\tdim(\chi)=1.$
By pulling the representation $\chi$ back
to a larger Galois group, we can assume 
$M$ is the ray class field modulo $\mathcal G$, where
$\mathcal G$ is principal and is divisible by the conductors of $\chi$ and the Hecke character $\psi$.  
Furthermore, since
$E$ is type (S), $\psi=\phi\circ \tnorm_{F/K}$ for some Hecke character $\phi$ of $K$.  

We begin by considering partial $L$-functions, which have only the Euler factors prime to $\mathcal G$.  By abuse of notation
the dependence on $\mathcal G$ is suppressed.
Then we have
\begin{multline*}
L(s,\overline{\psi}\otimes \chi)=
L(s,\overline{\phi}\circ \tnorm_{F/K}\otimes \chi)=\\L(s,\overline{\phi}\otimes\tind_G^\Gamma(\chi))=
\prod_{i=1}^{h} L(s,\overline{\phi}\otimes\chi_i)=\\
\prod_{i=1}^{h} \sum_{\gamma\in\Gamma}
\chi_i(\gamma)L(s,\overline{\phi},\gamma)
\end{multline*}
where $\tind_G^\Gamma(\chi)$ decom\-poses as $\oplus\chi_i$ and we have the par\-tial $L$-\-functions associated to an
Artin symbol $\gamma=[*,M/K]$ in $\Gamma$
\[
L(s,\overline{\phi},\gamma)=\sum_{\substack{\mathcal C\subset\mathcal O_K\\
[\mathcal C,M/K]=\gamma}} \overline{\phi}(\mathcal C)N(\mathcal C)^{-s}. 
\]
 This gives the $h$-th derivative at $s=0$ as
\begin{equation}\label{E:ab1}
L^{(h)}(0,\overline{\psi}\otimes \chi)=
\prod_{i=1}^{h} \sum_{\gamma\in\Gamma}\chi_i(\gamma) L^\prime(0,\overline{\phi},\gamma).
\end{equation}
In the next subsection we will give a generalization of what is usually called the Fro\-benius deter\-minant relation (act\-ually
due to Dedekind.)

\subsection{Generalized Dedekind determinant}

Suppose $\Gamma$ is a finite abel\-ian group, $G<\Gamma$, and
$
\chi:G\to\mathbb C^{\times}
$
is a character.  Write
\[
\pi=\tind_G^\Gamma(\chi)=\oplus \chi_i
\]
for the induced representation.  Fix once and for all a set $\boldsymbol{S}$ of coset representatives for
$G\backslash\Gamma$.   Let
$W$ be the vector space the induced representation acts in:
\[
W=\{F:\Gamma\to\mathbb C|F(\sigma x)=\chi(\sigma)F(x), \forall \sigma\in G,x\in\Gamma\},
\]
where $\Gamma$ acts by multiplication: $\pi(\gamma)F(x)=F(x\gamma)$.
Let 
\[
f:\Gamma \to\mathbb C
\]
 be any function, and define the operator on $W$
\[
\pi(f)=\sum_{\gamma\in\Gamma}f(\gamma)\pi(\gamma)
\]
The characters $\chi_i$ form a basis of $W$, and so do the characteristic functions of cosets $F_{G\gamma}$, where
\[
F_{G\gamma}(\sigma^\prime\gamma^\prime):=\begin{cases}
    \chi(\sigma^\prime)&\text{if } \gamma=\gamma^\prime\\
    0&\text{otherwise }
 \end{cases}
\]
The Dedekind determinant relation computes the determinant of $\pi(f)$ with respect to these two canonical
bases of $W$:
\begin{lemma}
\begin{align*}
\det\pi(f)&=\prod_{i=1}^{[\Gamma:G]}\sum_{\gamma \in\Gamma}f(\gamma)\chi_i(\gamma)\\
&=\det\left[\sum_{\tau\in G}\chi(\tau) f(\tau\gamma^\prime\gamma^{-1})\right]_{\gamma,\gamma^\prime\in
\boldsymbol{S}}
\end{align*}
\end{lemma}
\emph{Proof.}
The functions $\chi_i$ are eigenvectors of each $\pi(\gamma)$, with eigenvalue $\chi_i(\gamma)$, thus also
 eigenvectors of
$\pi(f)$ with eigenvalue 
\[
\sum_{\gamma\in\Gamma} f(\gamma)\chi_i(\gamma), 
\]
and so the first formula for the determinant is clear.    
Relative to our set of coset representatives $\boldsymbol{S}$, 
define a function (a factor set)
\begin{gather*}
\boldsymbol{S}\times \boldsymbol{S}\to G\\
(\gamma,\gamma^\prime)\mapsto g(\gamma,\gamma^\prime)
\end{gather*}
so that
\[
g(\gamma,\gamma^\prime)\gamma=\gamma^\prime \gamma^{\prime\prime}
\qquad\text{ where }\qquad
G\gamma=G\gamma^\prime G\gamma^{\prime\prime} 
\]
in the quotient $G\backslash\Gamma$.
An explicit computation shows that
\[
\pi(f)F_{G\gamma}=\sum_{\gamma^{\prime\prime}\in \boldsymbol{S}}\sum_{\sigma\in
G}\chi(\sigma g(\gamma,\gamma^\prime))f(\sigma \gamma^{\prime\prime})
F_{G\gamma^\prime},
\]
with $\gamma,\gamma^{\prime},\gamma^{\prime\prime}$ related as above.   A change of variables give the
matrix coefficients of $\pi(f)$ as in the lemma.
\qed
\begin{remark}
The case when $G$ is the trivial subgroup, and $\pi$ is the regular representation of $\Gamma$ is
the usual Frobenius determinant relation.
\end{remark}

Taking $L^\prime(0,\overline{\phi},\gamma)$ for the function $f(\gamma)$ in the Dedekind determinant, and using
(\ref{E:ab1}) realizes the $L$-function value as a determinant:
\begin{equation}\label{E:ab1more}
L^{(h)}(0,\overline{\psi}\otimes \chi)=
\det\left[
\sum_{\tau \in G}\chi(\tau)
L^\prime(0,\overline{\phi},\tau\gamma^\prime\gamma^{-1})
\right]_{\gamma,\gamma^\prime}.
\end{equation}
This application of the Dedekind determinant has long been a key ingredient for special values of $L$-functions.
For example, results on the conjecture of Birch and Swinnerton were obtained in
\cite{GS} and \cite{deS}.

\subsection{Partial $L$-functions and Kronecker series}

In this subsection we present, for completeness, a calculation of the derivative at $s=0$ of a partial $L$-function as
the value at $s=2$ of a Kronecker series.   

Since $\Phi_1(j(E))=j(\mathcal O_K)$, we have, for the lattice $\Lambda$ corresponding to $E_{\Phi_1}$
\[
\Lambda=\Omega\mathcal O_K
\]
for some $\Omega\in\mathbb C^\times$.  Further, for $\mathcal A$ an ideal of $\mathcal O_K$, we have $E^{[\mathcal
A,M/K]}$ is defined over $F$, corresponding to a lattice
\[
\Lambda_{\mathcal A}=h(\mathcal A)\Omega\mathcal A^{-1}
\]
for some $h(\mathcal A)$.  
The isogeny between $E$ and $E^{[\mathcal A,M/K]}$ is just multiplication by $h(\mathcal A)$.
Our hypothesis that $E$ is type (S) implies that $E$ is isogenous over $F$ to all its Galois conjugates,
which gives that $h(\mathcal A)\in F$ for all $\mathcal A$.  By composing isogenies one sees that $h$ is a crossed
homomorphism:
\[
h(\mathcal A^\prime\mathcal A)=h(\mathcal A^\prime)^{[\mathcal A,M/K]}h(\mathcal A).
\]
For more details on curves of type (S) (used throughout this section) see \cite{GS}.

Choose a set of representative $\mathcal A\in \boldsymbol{A}$ for the ideal class group of $K$, 
all prime to our fixed $\mathcal G$.
Choose also a fixed set of
representatives $\mathcal B\in\boldsymbol{B}$ so the Artin symbols $[\mathcal B, M/K]$ give every element of
$G=\tgal(M/F)$.  The ideals in $\boldsymbol{B}$ are principal as $F$ is the Hilbert class field, and we have that $\mathcal
B=(\phi(\mathcal B))$.  A given element  $\gamma\in\Gamma$ is then of the form $[\mathcal A\mathcal B,M/K]$ for some
$\mathcal A$ and some $\mathcal B$.  We get all ideals in this class by summing over $\alpha$ in $\mathcal
A^{-1}\mathcal G$ since then
\[
\phi(\mathcal B)+\alpha\equiv\phi(\mathcal B) \mod\mathcal G
\]
and
\[
[\mathcal A(\phi(\mathcal B)+\alpha),M/K]=[\mathcal A(\phi(\mathcal B)),M/K].
\]
So
\[
L(s,\overline{\phi},[\mathcal A\mathcal B,M/K])=\sum_{\alpha\in\mathcal A^{-1}\mathcal G}
\frac{\overline{\phi(\mathcal A)}}{N(\mathcal A)^s}
\frac{\overline{\phi((\phi(\mathcal B)+\alpha))}}{N(\phi(\mathcal B)+\alpha)^s}.
\]
Now in general we have 
\[
\phi((\lambda))=\phi_{\text{fin}}(\lambda)\lambda
\]
where $\phi_{\text{fin}}(\lambda)\in K$ only depends on $\lambda\mod\mathcal G$.  Since $\mathcal B=(\phi(\mathcal B))$
we get that $\phi_{\text{fin}}(\phi(\mathcal B))=1$, and
\[
\overline{\phi}((\phi(\mathcal B)+\alpha))=
\overline{\phi}_{\text{fin}}(\phi(\mathcal B)+\alpha)\overline{\phi(\mathcal B)+\alpha}=
\overline{\phi(\mathcal B)+\alpha}.
\]
Thus
\begin{equation}\label{E:ab2}
L(s,\overline{\phi},[\mathcal A\mathcal B,M/K])=\frac{\overline{\phi(\mathcal A)}}{N(\mathcal A)^s}
\sum_{\alpha\in\mathcal A^{-1}\mathcal G}
\frac{\overline{\phi(\mathcal B)+\alpha}}{|\phi(\mathcal B)+\alpha|^{2s}}.
\end{equation}

Let $\nu\in\Omega K^\times$ so that $(\nu/\Omega)=\mathcal G^{-1}$; i.e., $\nu$ is a $\mathcal G$ torsion point on
$\mathbb C/\Lambda$.  We see
that
\begin{equation}\label{E:ab3}
\frac{\overline{h(\mathcal A)\nu}}{|h(\mathcal A)\nu|^{2s}}
\sum_{\alpha\in\mathcal A^{-1}\mathcal G}
\frac{\overline{\phi(\mathcal B)+\alpha}}{|\phi(\mathcal B)+\alpha|^{2s}}=
K_1(\phi(\mathcal B)h(\mathcal A)\nu,0,s,\Lambda_{\mathcal A})
\end{equation}
since $\alpha\in\mathcal A^{-1}\mathcal G$ exactly when
$\omega=h(\mathcal A)\nu\alpha\in h(\mathcal A)\Omega\mathcal A^{-1}=\Lambda_{\mathcal A}.$
We combine equations (\ref{E:ab2}) and (\ref{E:ab3}),
multiply by $\Gamma(s)$ and use the functional equation (\ref{E:fneqn}) to see that
\begin{multline*}
\Gamma(s)L(s,\overline{\phi},[\mathcal A\mathcal B,M/K])=\\
\frac{\overline{\phi}(\mathcal A)}{N(\mathcal A)^s}
\frac{|h(\mathcal A)\nu|^{2s}}{\overline{h(\mathcal A)\nu}}
A(\Lambda_{\mathcal A})^{2-2s}
\Gamma(2-s)\times\\
K_1(0,\phi(\mathcal B)h(\mathcal A)\nu,2-s,\Lambda_{\mathcal A}).
\end{multline*}

Thus the partial $L$-function derivative at $s=0$ is computed by the Kronecker series value at $s=2$, which we are denoting
$K_{2,1}$:
\begin{lemma}\label{L:ab4}
\[
L^\prime(0,\overline{\phi},[\mathcal A\mathcal B,M/K])=
\frac{\overline{\phi}(\mathcal A)}{\overline{h(\mathcal A)\nu}}
A^2(\Lambda_{\mathcal A})K_{2,1}(\phi(\mathcal B)h(\mathcal A)\nu,\Lambda_{\mathcal A}).
\]
\end{lemma}

\subsection{}

For convenience we now number our  representatives $\mathcal A_i\in \boldsymbol{A}$ for the ideal class group of $K$, 
and choose them so  that 
$\mathcal A_i=\overline{\mathcal A_i}$ if $j(\mathcal A_i)$ is real ($1\le i\le r$), and
$\mathcal A_{i+s}=\overline{\mathcal A_i}$ if
$j(\mathcal A_i)$ is complex ($r<i\le s$).

Consider the matrix on the right hand side of equation \ref{E:ab1more}.  In terms of our representatives, we have
\[
\det\left[
\sum_{\mathcal B}\chi([\mathcal B,M/K])
L^\prime(0,\overline{\phi},[\mathcal A_i^{-1}\mathcal A_j \mathcal B,M/K])
\right]_{i,j}
\]

To use the relation between these partial $L$-function derivatives and Kronecker series above, we need to change  $\mathcal
A_i^{-1}$ to $\mathcal A_i.$ 
This induces a permutation of the ideal classes as well as an extra term from the representatives of the principal ideals in each
row. The permutation of the rows only
changes the determinant by $\pm1$, but the 
extra principal ideal needs to be absorbed from each row by a
change of variables in the sum, which alters the determinant.   This gives the $L$-function value as
\[
a(\chi)\det\left[
\sum_{\mathcal B}\chi([\mathcal B,M/K])
L^\prime(0,\overline{\phi},[\mathcal A_i\mathcal A_j \mathcal B,M/K])
\right]_{i,j}
\]
where 
$a(\chi)$ is the product of all terms introduced by these change of variables.
Note that, as the conjecture will require, 
\[
a(\chi)\in\mathbb Q(\chi),\qquad a(\chi^\sigma)=a(\chi)^\sigma\quad\forall
\sigma\in\tgal(\mathbb Q(\chi)/\mathbb Q).
\]
We use Lemma \ref{L:ab4} (with $\mathcal A_i\mathcal A_j$ instead of $\mathcal A$) on each entry of the matrix.  
Factor
$\overline{\phi}(\mathcal A_i)$ from row $i$ for each $i$, and similarly
$\overline{\phi}(\mathcal A_j)$ from each column.  Note that 
$\prod_{i}\overline{\phi}(\mathcal A_i)^2$ is in $K^\times$, since 
$\prod_{i}\mathcal A_i^2$ is principal.  
In fact in the \lq other half\rq of the $L$ function (coming from $\phi$ instead of
$\overline{\phi}$) we will see the complex conjugate of these terms.  So modulo $\mathbb Q^\times$ we can ignore them.
We have so far shown that
\begin{multline*}
L^{(h)}(0,\overline{\psi}\otimes \chi)\approx_{\mathbb Q(\chi)}\\
\det\left[
\sum_{\mathcal B}\chi([\mathcal B,M/K])
\frac{A^2(\Lambda_{\mathcal A_i\mathcal A_j})}{\overline{h(\mathcal A_i\mathcal A_j)\nu}}
K_{2,1}(\phi(\mathcal B)h(\mathcal A_i\mathcal A_j)\nu,\Lambda_{\mathcal A_i\mathcal A_j})
\right]_{i,j}
\end{multline*}

We can now apply the distribution relation of Proposition \ref{P:distribution} to the isogeny 
\[
h(\mathcal A_i)^{[\mathcal A_j, F/K]}:E^{[\mathcal A_j, F/K]}\to E^{[\mathcal A_i\mathcal A_j, F/K]}
\]
of degree $N(\mathcal A_i)$.
We need to use the fact that in general
\[
A(\Lambda_{\mathcal A})=A(h(\mathcal A)\Omega\mathcal A^{-1})=
\frac{|h(\mathcal A)\Omega|^2\sqrt{|D|}}{2\pi N(\mathcal A)}
\]
where $D$ is the discriminant of $K$.  So
\[
A^2(\Lambda_{\mathcal A_i\mathcal A_j})=\frac{|h(\mathcal A_i)^{[\mathcal A_j, F/K]}|^4}{N(\mathcal
A_i)^2}A^2(\Lambda_{\mathcal A_j}).
\]
In summary, we've got
\begin{multline*}
L^{(h)}(0,\overline{\psi}\otimes \chi)\approx_{\mathbb Q(\chi)}\\
\det\left[
\sum_
{\substack{\mathcal B\in\boldsymbol{B}\\t\in\tker h(\mathcal A_i)}}
\chi([\mathcal B,M/K])
\frac{A^2(\Lambda_{\mathcal A_j})}{\overline{h(\mathcal A_j)\nu}}
K_{2,1}(\phi(\mathcal B)h(\mathcal A_j)\nu-t,\Lambda_{\mathcal A_j})
\right]_{i,j}
\end{multline*}
By class field theory, the elements $\mathcal A_j$ of the class group correspond via our fixed embedding $\Phi_1$ to the other
embeddings $\Phi_j$ which have the same restriction to $K$.  That is
\[
\mathbb C/\Lambda_{\mathcal A_j}=E^{[\mathcal A_j,M/K]}=E_{\Phi_j}.
\]
We use the homothety property to factor a scalar $h(\mathcal A_j)\Omega$ out of 
the $K_{2,1}$ in each column,  to convert $\Lambda_{\mathcal A_j}$ to $\mathcal A^{-1}_j$.  If we take our representatives
to be integral ideals not divisible by any rational integer but 1, it is easy to see that
$\mathcal A^{-1}_j$ has a lattice basis of the form $[1,\tau_j]$ as required.  

The divisor associated to the sum over the torsion points
in the kernel of $h(\mathcal A_i)$ comes from a symbol $\xi_i$ in $\mathbb Q K_2(E_M)$
by the theorem of Bloch \cite{Birv}.  Multiplication by $\phi(\mathcal B)$ on the torsion points of one of these curves acts by
Galois automorphism
$[\mathcal B,M/K]$.  Recall these fix $F$, and since our isogenies $h(\mathcal A_i)$ are defined over $F$ we can re-write
the sum over $\mathcal B\in\boldsymbol{B}$ as a sum over $\tau\in G$.  Thus we see
\[
L^{(h)}(0,\overline{\psi}\otimes \chi)\approx_{\mathbb Q(\chi)}
\det\left[\overline{\Omega/\nu}
\sum_
{\tau\in G}
\chi(\tau)
\treg(\xi_i^\tau)_{\Phi_j}
\right]_{i,j}
\]

Applying this construction to $\psi$ instead of $\overline{\psi}$, we get the same formula but with the conjugate embeddings
$\overline{\Phi}_j$ instead of $\Phi_j$.  However,
we have so far only constructed $h$ symbols in $K_2$ which we can relate to the $L$-value.  To get $2h$ symbols, we take
advantage of the fact that $K_2(E_M)$ is an $\mathcal O_K$ module.  More specifically, the point is that we can vary the ideal
$\mathcal G$, changing the matrix of partial $L$ functions and thus also the matrix of symbols, by a constant in $K$.  
For any ideal $\mathcal P$, one
sees that  for any $\gamma$ in $\Gamma$, the partial $L$-functions satisfy
\[
L_{\mathcal G}(s,\overline{\phi},\gamma)=
L_{\mathcal G\mathcal P}(s,\overline{\phi},\gamma)+
\overline{\phi}(\mathcal P)N(\mathcal P)^{-s}
L_{\mathcal G}(s,\overline{\phi},\gamma\cdot[\mathcal P,M/K]^{-1})
\]
where we now, of course, need to keep track of the ideal in the notation for the partial $L$-function.
Thus if $\mathcal P$ is principal
\begin{multline*}
\sum_{\tau\in G}\chi(\tau)L_{\mathcal G\mathcal P}^\prime(0,\overline{\phi},\gamma\tau)=\\
(1-\overline{\phi}(\mathcal P)\chi([\mathcal P,M/K]))\times
\sum_{\tau\in G}\chi(\tau)L_{\mathcal G}^\prime(0,\overline{\phi},\gamma\tau).
\end{multline*}
Choose
two principal prime ideals $\mathcal P^+$ and $\mathcal P^-$, and let $\mathcal G^{\pm}=\mathcal G\mathcal P^{\pm}$.  Let
\[
\pi^{\pm}=1-\overline{\phi}(\mathcal P^{\pm})\chi([\mathcal P^{\pm},M/K]).
\]
When we change $\mathcal G$ to $\mathcal G^+$ or $\mathcal G^-$,  the matrix on the right 
hand side of (\ref{E:ab1more}) changes by the scalar $\pi^+$ or $\pi^-$.  Following the calculation of the previous section to the
end, we see the same is true for the matrix of regulators of symbols supported on the torsion.  That is, let
\[
R=\left[
\sum_
{\tau\in G}
\chi(\tau)
\treg(\xi_i^\tau)_{\Phi_j}
\right]_{i,j}
\]
the matrix corresponding to symbols on the $\mathcal G$ torsion, then
\[
\pi^{\pm}R=\left[
\sum_
{\tau\in G}
\chi(\tau)
\treg(\xi_{i,\pm}^\tau)_{\Phi_j}
\right]_{i,j}
\]
where the symbols $\xi_{i,\pm}$ come from the $\mathcal G^{\pm}$ torsion.  We have shown above that
\[
\det(R)\det(\overline{R})\approx_{\mathbb Q(\chi)}
L^{(h)}(0,\overline{\psi}\otimes \chi)
L^{(h)}(0,\psi\otimes \chi)=
L^{(2h)}(0,E\otimes\chi).
\]
By Proposition \ref{P:case1reg}, the elliptic Stark conjecture will follow from the linear algebra
\begin{lemma} For an $h\times h$ matrix $R$, and scalars $\pi^+,\pi^-$ in $\mathbb Q(\sqrt{D})$, 
\[
\det\begin{bmatrix}
\pi^+R&\overline{\pi^+R}\\
\pi^-R&\overline{\pi^-R}
\end{bmatrix}=\kappa \det(R)\det(\overline{R}),\quad
\kappa\in
	\begin{cases}
    \mathbb Q&\text{ if }h\text{ is even }\\
    \mathbb Q\cdot\sqrt{D}&\text{ if }h\text{ is odd }
 \end{cases}
\]
\end{lemma}
\emph{Proof.}   This is the fancy version of the Laplace expansion theorem, where we are expanding on the first $h$ columns. 
Note that in the sum over $h$ by $h$ matrices, we include row $i$ of $\pi^+R$ if and only if we omit row $i$ of
$\pi^-R$, as these are the only terms with nonzero determinant.  This gives $\det(R)\det(\overline{R})$ times a sum of powers
of $\pi^+$ and $\pi^-$.  If $h$ is even, the Laplace theorem gives that each term occurs with sign $+1$, and it is a trace from
$\mathbb Q(\sqrt{D})$ to $\mathbb Q$.  If $h$ is odd, then group complementary terms in the sum, which necessarily occur
with opposite sign.  One sees that it is a \lq skew trace \rq, so is in $\mathbb Q\cdot\sqrt{D}$.\qed

\subsection{The simple zero}

In the Stark conjectures, the case when the $L$-function has a simple zero gets special attention.  We observe here
that if $[F:\mathbb
Q]\cdot\tdim(V)=1$, then $F=\mathbb Q$.  Since $E$ is defined over $F$, we must necessarily have the class
number $h(K)=1$.  Since $\tdim(V)=1$, $\chi$ factors through an abel\-ian extension $\tgal(M/\mathbb Q)$,
and $M$ contains $K$ by hypothesis.  By the remark at the beginning of \S \ref{S:Torsion}, $E_L$ is of type
(S) for any intermediate field $L$ with $\mathbb Q\subseteq L\subseteq M$.  Assuming the Dimension conjecture, the results in
the previous subsection give the elliptic Stark conjecture in this case.

\section{Curves not type (S)}\label{S:NotS}
The positive results so far have all been for curves of type (S), even though we introduced the elliptic Stark conjecture to
study the general case.  In this section we remedy this defect with the following
\begin{theorem}
If $F$ is abelian over $K$, and $E$ is any elliptic curve over $F$ with
complex multiplication by $\mathcal O_K$, then there is a Galois extension $M$ of $F$, and a character $\chi$ of
$\tgal(M/F)$ such that the elliptic Stark conjecture holds for $L(s,E\otimes\chi).$
\end{theorem}
\begin{remark}
This theorem does not assume the \lq Dimension conjecture\rq of \S\ref{S:Intro}.
\end{remark}
\emph{Proof}
By
\cite{Robert} Corollaire 2 we know there exists an elliptic curve $E^\prime$ defined over $F$ which
is of type (S).   By Theorem 9.1.3 of \cite{Gr}, we may assume that
$E^\prime$ and
$E$ have the same $j$ invariant.  Thus $E^\prime$ is a model of $E$ and we can write Weierstrass equations
\begin{align*}
E:\quad y^2=&4x^3-g_2x-g_3\\
E^\prime:\quad y^2=&4x^3-d^2g_2x-d^3g_3
\end{align*}
with $d$ in $F$.  The curves $E$ and $E^\prime$ become isomorphic over $M=F(\sqrt{d})$ via
\begin{align*}
\phi:E&\to E^\prime\\
(x,y)&\mapsto (x^\prime=dx,y^\prime=d^{3/2}y)
\end{align*}
We get a map on functions fields and a map on $K$-groups
\[
\phi^*:\mathbb QK_2(E^\prime_F)\to\mathbb QK_2(E_F),
\]
as in \S\ref{S:Analysis}.  

The inclusions $i:F(E)\to M(E)$ and $i^\prime:F(E^\prime)\to M(E)$ also give maps 
\begin{gather*}
i^*:\mathbb QK_2(E_F)\to\mathbb QK_2(E_M)\\
i^{\prime*}:\mathbb QK_2(E^\prime_F)\to \mathbb QK_2(E_M)
\end{gather*} 

Note that the triangle formed by these three
maps does not commute.  In fact $M(E)=M(E^\prime)$ is a Galois extension of $F(x)=F(x^\prime)$, with Galois group the Klein
four group.  The quadratic subfields are $M(x)$, $F(E)$, and $F(E^\prime)$.  The subgroups of order two which fix these fields
are generated by
$[-1]$, $\tau$, and $\tau^\prime$, where
\[
\tau:
	\begin{cases}
    x&\to x\\
    y&\to y\\
\sqrt{d}&\to-\sqrt{d}
 \end{cases}\qquad
\tau^\prime:
	\begin{cases}
    x^\prime&\to x^\prime\\
    y^\prime&\to y^\prime\\
\sqrt{d}&\to-\sqrt{d}
 \end{cases}
\]
Of course, $\tau^\prime=\tau\circ[-1]$, and both $\tau$ and $\tau^\prime$ restrict to the same nontrivial automorphism of $M$
over $F$.  For $f$ in $F(E^\prime)$, we see 
\[
(\phi^*f)^\tau=\phi^*f,\qquad f^\tau=f\circ[-1]=[-1]^*f.  
\]
Thus for symbols $\xi$ in
$\mathbb QK_2(E^\prime_F)$, 
\[
(i^*\phi^*\xi)^\tau=i^*\phi^*\xi,\qquad (i^{\prime*}\xi)^\tau=[-1]^*i^{\prime*}\xi.
\]

Since $E^\prime$ is type (S) there exist $n=[F:\mathbb Q]$ symbols $\xi_i$ such that the
\lq $L$-value conjecture\rq of \S\ref{S:Intro} is true.   In the notation of \S\ref{S:Stark},
this says that
\[A(E^\prime_F)=R(E^\prime_F)/c(E^\prime_F)\] is rational.
Let $\chi$ be the nontrivial character of $\tgal(M/F)$.  Since 
\[
L(s,E_F)L(s,E^\prime_F)=L(s,E_M)=L(s,E_F)L(s,E_F\otimes\chi)
\]
we see that $c(E^\prime_F)=c(E_F,\chi).$  

It remains to relate 
$R(E^\prime_F)$ to $R(E_F,\chi).$  For notational convenience we suppress the 
inclusions $i^*$ and $i^{\prime*}$.  We see that for symbols $\xi_i$ as above,
\[
(\frac{\phi^*\xi_i+\xi_i}{2})^\tau=\frac{\phi^*\xi_i-\xi_i}{2},
\]
so the $n$ symbols $(\phi^*\xi_i+\xi_i)/2$ satisfy the requirement of \S\ref{S:Stark} to be a set 
$\boldsymbol{M}$ of \lq Minkowski symbols\rq.  
Fix any embedding $\Phi_j$ of $F$ into $\mathbb C$.  In the formula for $R(E_F,\chi,f_{\boldsymbol{M}})$ in Proposition
\ref{P:case1reg}, the $i,j$ entry in the sum over the Galois group simplifies as
\[
\treg(\frac{\phi^*\xi_i+\xi_i}{2})_{\Phi_j}-\treg((\frac{\phi^*\xi_i+\xi_i}{2})^\tau)_{\Phi_j}=\treg(\xi_i)_{\Phi_j}.
\]
Thus $R(E^\prime_F)=R(E,\chi)$ and so 
\[
A(E_F,\chi)=R(E_F,\chi)/c(E_F,\chi)=R(E^\prime_F)/c(E^\prime_F)
\]
is rational.
\qed

\end{document}